\documentclass{amsart}
\title{Topological games and covering dimension}
\author{Liljana Babinkostova}
\usepackage{amsfonts}
\newcommand{\sone}{{\sf S}_1}
\newcommand{\gone}{{\sf G}_1}

\newcommand{\Sc}{{\sf S}_c}
\newcommand{\Gc}{{\sf G}_{c}}
\newcommand{\op}{\mathcal{O}}
\newcommand{\opkfd}{{\mathcal{O}}_{\sf kfd}}
\newcommand{\om}{\Omega}

\newcommand{\opfd}{\mathcal{O}_{\sf fd}}

\newcommand{\naturals}{{\mathbb N}}

\newcommand{\smirnovspace}{{\sf S}_{\omega}}
\newcommand{\hilbertcube}{{\mathbb H}}
\newcommand{\polspace}{{\mathbb K}}
\newcommand{\closedunitinterval}{{\mathbb I}}
\newtheorem{theorem}{{\bf Theorem }}
\newtheorem{lemma}[theorem]{{\bf Lemma }}
\newtheorem{corollary}[theorem]{{\bf Corollary}}
\newtheorem{proposition}[theorem]{{\bf Proposition}}	
\newtheorem{problem}{{\bf Problem }}

\newtheorem{conjecture}{{\bf Conjecture}}

\begin{document}
\maketitle
\begin{abstract}
We consider a natural way of extending the Lebesgue covering dimension to various classes of infinite dimensional separable metric spaces.   
\end{abstract}

All spaces in this paper are assumed to be separable metric spaces. Infinite games can be used in a natural way to define ordinal valued functions on the class of separable metric spaces. One of our examples of such an ordinal function coincides in any finite dimensional metric spaces with the covering dimension of the space, and may thus be thought of as an extension of Lebesgue covering dimension to all separable metric spaces. We will call this particular extension of Lebesgue covering dimension the \emph{game dimension} of a space.

Game dimension is defined using a game motivated by a selection principle. Several natural classical selection principles are related to the one motivating game dimension, and their associated games can be used to compute upper bounds on game dimension. These games, and the upper bounds they provide, are interesting in their own right. We develop two such games and use them then to obtain upper bounds on our game dimension.

We also compute the game dimension of a few specific examples. 

\section{Selection principles, open covers and games}

The selection principle $\sone(\mathcal{A},\mathcal{B})$ states that there is for each sequence $(A_n:n\in\naturals)$ of elements of $\mathcal{A}$ a sequence $(B_n:n\in\naturals)$ such that for each $n$ we have $B_n\in A_n$, and $\{B_n:n\in\naturals\}\in\mathcal{B}$.

The selection principle $\Sc(\mathcal{A},\mathcal{B})$ states that there is for each sequence $(A_n:n\in\naturals)$ of members of $\mathcal{A}$ a sequence $(B_n:n\in\naturals)$ of sets such that for each $n$, $B_n$ is a pairwise disjoint family of sets and refines $A_n$ and $\bigcup_{n\in\naturals}B_n\in\mathcal{B}$.

For a collection $\mathcal{T}$ of subsets of a topological space $X$ we call an open cover $\mathcal{U}$ of $X$ a ``$\mathcal{T}$-cover" if for each $T\in \mathcal{T}$, there is a $U\in\mathcal{U}$ with $T\subseteq U$. The symbol $\op(\mathcal{T})$ denotes the collection of $\mathcal{T}$-covers of $X$. A trivial situation, but one we cannot ignore, arises when $X$ itself is a member of $\mathcal{T}$. We don't follow the usual practice of requiring $X\not\in\mathcal{U}$ for $\mathcal{U}\in\op(\mathcal{T})$. The motivation, as will be seen below, is that allowing this trivial situation provides a uniformity to the statements of some of our results.

A few special combinatorial properties of the family $\mathcal{T}$ are important for our considerations. Here are some of them: $\mathcal{T}$ is said to be 
\emph{up-directed} if for all $A$ and $B$ in $\mathcal{T}$ there is a $C$ in $\mathcal{T}$ with $A\bigcup B \subset C$.
It is said to be \emph{first-countable} if there is for each $T\in\mathcal{T}$ a sequence $(U_n:n\in\naturals)$ of open sets such that for each $n$, $U_n\supset U_{n+1}\supset T$, and for each open set $V\supset T$, there is an $n$ with $V\supseteq U_n$.
We shall say that $X$ is $\mathcal{T}$-\emph{first countable} if   there is for each $T\in\mathcal{T}$ a sequence $(U_n:n=1,2,\ldots)$ of open sets such that for all $n$, $T\subset U_{n+1}\subset U_n$,   and for each open set $U\supset T$ there is an $n$ with $U_n\subset U$. Let $\langle\mathcal{T}\rangle$ denote the subspaces which are unions of countably many elements of $\mathcal{T}$.
When $\mathcal{T}$ is a collection of compact sets in a metrizable space $X$, $\mathcal{T}$ is first countable and up-directed. Call a subset $\mathcal{C}$ of $\mathcal{T}$ \emph{cofinal} if there is for each $T\in\mathcal{T}$ a $C\in\mathcal{C}$ with $T\subseteq C$. 

When $\mathcal{T}$ is the collection of one-element subsets of $X$, then $\op(\mathcal{T})$ will be denoted simply $\op$. With $\opfd$ we denote the collection of all finite dimensional subsets of a separable metric space and with $\opkfd$ we denote the collection of all compact, finite dimensional subsets of a separable metric space. When $\mathcal{T}$ is the collection of finite subsets of $X$ then $\mathcal{T}$-covers are called $\omega$-covers, and $\om$ denotes $\op(\mathcal{T})$.

The collection of open covers having just two elements is denoted $\op_2$. A topological space is said to be weakly infinite dimensional if it has the property $\Sc(\op_2,\op)$. The class of spaces satisfying $\Sc(\op,\op)$ was introduced in \cite{addisgresham}. If a space is a union of countably many zerodimensional subsets it is said to be countable dimensional. Hurewicz introduced the latter notion. As was shown in \cite{addisgresham}, countable dimensional implies $\Sc(\op,\op)$, and $\Sc(\op,\op)$ implies $\Sc(\op_2,\op)$. If a space is not weakly infinite dimensional, then it is said to be \emph{strongly infinite dimensional}. The Hilbert cube $\hilbertcube$ is an example of strongly infinite dimensional space.

Let $\alpha$ be an ordinal number. The game $\gone^{\alpha}(\mathcal{A},\mathcal{B})$ associated with the selection principle $\sone(\mathcal{A},\mathcal{B})$ is as follows: The players play an inning per $\gamma<\alpha$. In the $\gamma$-th inning ONE first chooses an $A_{\gamma}\in\mathcal{A}$: TWO then responds with a $B_{\gamma}\in A_{\gamma}$. A play $A_0,\, B_0,\, \cdots,\, A_{\gamma},\, B_{\gamma},\, \cdots$ of length $\alpha$ is won by TWO if $\{B_{\gamma}:\gamma<\alpha\}\in\mathcal{B}$. Else, ONE wins.

When for a set S and families $\mathcal{A}$ and $\mathcal{B}$ there is an ordinal number $\alpha$ such that TWO has a winning strategy in the game $\gone^{\alpha}(\mathcal{A},\mathcal{B})$ played on S, then we define:
\[
  {\sf tp}_{\sone(\mathcal{A},\mathcal{B})}(S) = \min\{\alpha:\mbox{ TWO has a winning strategy in }\gone^{\alpha}(\mathcal{A},\mathcal{B}) \mbox{ on }S\}.
\]
We adopt the convention that ${\sf tp}_{\sone(\mathcal{A},\mathcal{B})}(\emptyset)=-1$.

The ordinal-valued function ${\sf tp}_{\sone(\mathcal{A},\mathcal{B})}(S)$ has been studied for a few specific choices of the families $\mathcal{A}$ and $\mathcal{B}$. Examples closely affiliated with what we will examine here can be found in \cite{lbms}, \cite{baldwin}, \cite{bernerjuhasz}, \cite{danielsgruenhage}, \cite{diaggamelength} and \cite{succtype}.

 The following monotonicity properties are easy to check:
\begin{lemma}\label{opvsom} For each separable metric space $X$, for each closed set $Y\subseteq X$ and for families $\mathcal{S}\subset\mathcal{T}$ and $\mathcal{A}\subset\mathcal{B}$ the following hold:
\begin{enumerate}
\item{${\sf tp}_{\sone(\op(\mathcal{T}),{\mathcal{A}})}(X) \le {\sf tp}_{\sone(\op(\mathcal{S}),{\mathcal{A}})}(X)$.}
\item{${\sf tp}_{\sone(\op(\mathcal{T}),{\mathcal{B}})}(X) \le {\sf tp}_{\sone(\op(\mathcal{T}),{\mathcal{A}})}(X)$.}
\item{${\sf tp}_{\sone(\op(\mathcal{T}),{\mathcal{A}})}(Y) \le {\sf tp}_{\sone(\op(\mathcal{T}),{\mathcal{A}})}(X)$.}
\end{enumerate}
\end{lemma}

In particular we find that ${\sf tp}_{\sone(\opfd,{\mathcal{A}})}(X) \le {\sf tp}_{\sone(\opkfd,{\mathcal{A}})}(X)$.

Some of the results from our investigation in \cite{lbms} have some applications in this paper. We first recall these.

In order to compare the ordinals ${\sf tp}_{S_{\xi}(\mathcal{A},\mathcal{B})}(X)$ for various choices $\xi$, $\mathcal{A}$ and $\mathcal{B}$, we extend Lemma 1 of \cite{lbms} as follows:

\begin{lemma}\label{setandopenset} Let $\alpha$ be an ordinal number, and let $F$ be a strategy for TWO in the game $\gone^{\alpha}(\op(\mathcal{T}),\op)$. Then there is for each $\nu<\alpha$, for each sequence $(O_{\beta}:\beta<\nu)$ from $\op(\mathcal{T})$ a set $C\in\mathcal{T}$ such that for each open set $U\supseteq C$ there is an $O\in\op(\mathcal{T})$ with $F((O_{\beta}:\beta<\nu)\frown(O)) = U$.
\end{lemma} 

\begin{theorem}[\cite{lbms} Theorem 4]\label{twoskerneldirected} Let $\alpha$ be an infinite ordinal and let $\mathcal{T}$ be up-directed. If $F$ is any strategy for TWO in $\gone^{\alpha}(\mathcal{O}(\mathcal{T}),\mathcal{O})$ and if $X$ is $\mathcal{T}$-first countable, then there is for each set $T\in\langle\mathcal{T}\rangle$ a set $S\in\langle\mathcal{T}\rangle$ such that:
$T\subseteq S$ and for any closed set $C\subset X\setminus S$, there is an $\omega$-length $F$-play 
\[
  O_1,\, T_1,\, \cdots,\, O_n,\, T_n\,\cdots
\]
such that $T\subseteq \bigcup_{n=1}^{\infty}T_n \subseteq X\setminus C$.
\end{theorem}

\begin{theorem}[\cite{lbms} Theorem 5]\label{twoskerneldirectedcofinal} Let $\alpha$ be an infinite ordinal and let $\mathcal{T}$ be up-directed. If $F$ is any strategy for TWO in $\gone^{\alpha}(\mathcal{O}(\mathcal{T}),\mathcal{O})$ and if $X$ is $\mathcal{C}$-first countable where $\mathcal{C}\subset\mathcal{T}$ is cofinal in $\mathcal{T}$, then there is for each set $T\in\langle\mathcal{T}\rangle$ a set $S\in\langle\mathcal{C}\rangle$ such that:
$T\subseteq S$ and for any closed set $C\subset X\setminus S$, there is an $\omega$-length $F$-play 
\[
  O_1,\, T_1,\, \cdots,\, O_n,\, T_n\,\cdots
\]
such that $T\subseteq \bigcup_{n=1}^{\infty}T_n \subseteq X\setminus C$.
\end{theorem}

\begin{lemma}[\cite{lbms} Theorem 6]\label{opomegalength} If $\mathcal{T}$ has a cofinal subset consisting of ${\sf G}_{\delta}$-sets, then TWO has a winning strategy in $\gone^{\omega}(\op(\mathcal{T}),\op)$ if, and only if, $X$ is a union of countably many elements of $\mathcal{T}$.
\end{lemma}

\begin{lemma}[\cite{lbms} Lemma 13]\label{twoskernelop} Let $X$ and $\mathcal{T}$ be such that $X\not\in\langle\mathcal{T}\rangle$ and $\mathcal{T}$ is up-directed and first-countable. If $\alpha$ is the least ordinal for which there is a $B\in\langle\mathcal{T}\rangle$ such that each closed set $C\subset X\setminus B$ satisfies ${\sf tp}_{\sone(\op(\mathcal{T}),\op)}(C)\le\alpha$, then ${\sf tp}_{\sone(\op(\mathcal{T}),\op)}(X)= \omega+\alpha$.
\end{lemma}
The converse is also true:
\begin{lemma}\label{structure} Let $X$ and $\mathcal{T}$ be such that $X\not\in\langle\mathcal{T}\rangle$ and $\mathcal{T}$ is up-directed and first-countable. If ${\sf tp}_{\sone(\op(\mathcal{T}),\op)}(X)= \omega+\alpha$ then there is a $B\in\langle\mathcal{T}\rangle$ such that each closed set $C\subset X\setminus B$ satisfies ${\sf tp}_{\sone(\op(\mathcal{T}),\op)}(C)\le\alpha$, and $\alpha$ is the minimal such ordinal.
\end{lemma}
{\flushleft{\bf Proof:}} For let $F$ be a winning strategy for TWO in the game $\gone^{\omega+\alpha}(\op(\mathcal{T}),\op)$. By Lemma \ref{setandopenset} choose $C_{\emptyset}\in\mathcal{T}$ such that for each open set $U\supseteq C_{\emptyset}$ there is a $\mathcal{U}\in\op(\mathcal{T})$ with $U = F(\mathcal{U})$. Since $\mathcal{T}$ is first countable choose a sequence $(U_n:n<\omega)$ of open sets with: For each open set $U\supseteq C_{\emptyset}$ there is an $n$ with $U\supseteq U_{n}\supseteq C_{\emptyset}$, and for all $n$, $U_n\supset U_{n+1}$. For each $n$ choose $O_n\in\op(\mathcal{T})$ with $U_n = F(O_n)$. Applying Lemma \ref{setandopenset} to each $O_n$, fix sets $C_n\in\mathcal{T}$ such that there is for each open $U\supseteq C_n$ an $O\in\op(\mathcal{T})$ with $U = F(O_n,O)$. Since $\mathcal{T}$ is first countable choose for each $n$ a sequence $(U_{n,m}:m<\omega)$ of open sets such that: For each open set $U\supseteq C_{n}$ there is an $m$ with $U\supseteq U_{n,m}\supseteq C_{n}$, and for all $m$, $U_{n,m}\supset U_{n, m+1}$. Then choose for each $m$ an $O_{n,m}\in\op(\mathcal{T})$ with $U_{n,m}=F(O_n,O_{n,m})$.

Continuing like this we find for each finite sequence $(n_0,\cdots,n_j)$ sets $C_{n_0,\cdots,n_j}\in\mathcal{T}$, and open sets $U_{n_0,\cdots,n_j,t}$, and $O_{n_0,\cdots,n_j,t}\in\op(\mathcal{T})$, $t<\omega$, such that for each open set $U\supseteq C_{n_0,\cdots,n_j}$ there is a $t$ with $U\supseteq U_{n_0,\cdots,n_j,t}\supseteq C_{n_0,\cdots,n_j}$ and $U_{n_0,\cdots,n_j,t} = F(O_{n_0},O_{n_0,n_1},\cdots,O_{n_0,\cdots,n_j,t})$, and for all $t$, $U_{n_0,\cdots,n_j,t}\supset U_{n_0,\cdots,n_j,t+1}$.

Put $B=\cup_{\tau\in\,^{<\omega}\omega}C_{\tau}$, an element of $\langle\mathcal{T}\rangle$. Consider any closed set $C\subset X\setminus B$. Since $U=X\setminus C\supseteq B$ is open, choose $n_0$ with $C_{\emptyset}\subseteq U_{n_0}\subset U$. Then choose $n_1$ with $C_{n_0}\subseteq U_{n_0,n_1}\subseteq U$, and then $n_2$ with $C_{n_0,n_1}\subseteq U_{n_0,n_1,n_2}\subseteq U$, and so on. In this way we obtain an $\omega$-sequence of moves during which TWO used the strategy $F$, and the union of TWO's responses, say $V$, is a subset of $U$. Since $F$ is a winning strategy for TWO in the game $\gone^{\omega+\alpha}(\op(\mathcal{T}),\op)$ it follows that ${\sf tp}_{\sone(\op(\mathcal{T},\op)}(C) \le {\sf tp}_{\sone(\op(\mathcal{T},\op)}(X\setminus U)\le {\sf tp}_{\sone(\op(\mathcal{T},\op)}(X\setminus V) = \alpha$. The closed set $X\setminus V$ witnesses that $\alpha$ is the minimal such ordinal. $\square$

\section{The ordinal ${\sf tp}_{\sone(\opkfd,\op)}(X)$.}

The properties of open point-covers of spaces are at the basis of Lebesgue's notion of covering dimension. The points have, in that theory, dimension zero. We will consider a more general situation of open $\mathcal{T}$-covers of spaces. In this case by analogy with point covers, we require that the notion of dimension assign the value 0 to the members of $\mathcal{T}$.  We feature here the specific case when $\mathcal{T}$ is the collection of finite dimensional compact spaces. 

{\flushleft 2.1 {\bf When ${\sf tp}_{\sone(\opkfd,\op)}(X)$ is countable}}\\

\begin{proposition}\label{omegalength} 
For X a separable metric space the following are equivalent:
\begin{enumerate}
\item{X is a countable union of compact finite dimensional subsets.}
\item{$1 < {\sf tp}_{\sone(\opkfd,\op)}(X) \le \omega$.}

\end{enumerate}
\end{proposition}
We shall see that (2) of Proposition \ref{omegalength} cannot be improved. \\

{\flushleft 2.2 {\bf Products and unions.}}\\

\begin{lemma}\label{cptfindimfactor} For $X$ a metric space, and $Y$ a compact finite dimensional metric space, ${\sf tp}_{\sone(\opkfd,\op)}(X\times Y) = {\sf tp}_{\sone(\opkfd,\op)}(X)$.
\end{lemma}
{\bf Proof:}First note that since $X$ is a closed subset of $X\times Y$, we have ${\sf tp}_{\sone(\opkfd,\op)}(X)\le {\sf tp}_{\sone(\opkfd,\op)}(X\times Y)$. We show that ${\sf tp}_{\sone(\opkfd,\op)}(X\times Y) \le {\sf tp}_{\sone(\opkfd,\op)}(X)$.

Let $F$ be a winning strategy for TWO in $\gone^{\alpha}(\opkfd,\op)$ played on $X$. TWO uses $F$ as follows to play the game on $X\times Y$. First observe that if $O$ is in $\opkfd$ then there is an $S(O)\in\opkfd$ such that each element of $S(O)$ is of the form $U\times Y$, and such that $S(O)$ is a refinement of $O$.

Define a strategy $G$ for TWO in the game $\gone^{\alpha}(\opkfd,\op)$ on $X\times Y$ as follows:
In inning 1, when ONE plays $O_1$, put $B_1 = \{U\subset X:U\times Y\in S(O_1)\}$, a member of $\opkfd$ for $X$. Then compute 
$W_1 = F(B_1)\in B_1$, and choose $T_1\in O_1$ with $W_1\times Y\subseteq T_1$. Define $G(O_1) = T_1$.

Let $\gamma<\alpha$ be given, as well as a sequence $(O_{\xi}:\xi\le \gamma)$ of members of $\opkfd$ for $X\times Y$. For each $\xi\le\gamma$ define $B_{\xi} = \{U\subset X:U\times Y\in S(O_{\xi})\}$, a member of $\opkfd$ on $X$. Compute $W_{\gamma} = F(B_{\xi}:\xi\le \gamma)\in B_{\gamma}$ and then choose $T_{\gamma}\in O_{\gamma}$ with $W_{\gamma}\times Y\subseteq T_{\gamma}$, and put $G(O_{\xi}:\xi\le\gamma)=T_{\gamma}$.

Then $G$ is a winning strategy for TWO in $\gone^{\alpha}(\opkfd,\op)$ on $X\times Y$. $\square$

\begin{lemma}\label{twowinsinfinite}
Let $X$ be a metric space which is not compact and finite dimensional. If ${\sf tp}_{\sone(\opkfd,\op)}(X)$ is finite then there is a nonempty finite dimensional compact set $C$ such that for each open set $U\supset C$ we have ${\sf tp}_{\sone(\opkfd,\op)}(X\setminus U) \le {\sf tp}_{\sone(\opkfd,\op)}(X)-1$, and for some open set $U\supset C$ we have ${\sf tp}_{\sone(\opkfd,\op)}(X\setminus U) = {\sf tp}_{\sone(\opkfd,\op)}(X)-1$     
\end{lemma}
{\flushleft{\bf Proof:}} Suppose not. For each compact finite dimensional $C\subset X$ choose an open set $U(C)$ such that $C\subset U(C)$, and ${\sf tp}_{\sone(\opkfd,\op)}(X\setminus U(C))\ge {\sf tp}_{\sone(\opkfd,\op)}(X)$. Then $O_1=\{U(C):C\subset X \mbox{ finite dimensional compact}\}\in\opkfd$, and for each strategy $G$ of TWO, ${\sf tp}_{\sone(\opkfd,\op)}(X\setminus G(O_1)) \ge  {\sf tp}_{\sone(\opkfd,\op)}(X)$. But this contradicts the fact that ${\sf tp}_{\sone(\opkfd,\op)}(X)$ is finite.

Suppose that for each compact finite dimensional $C\subset X$, for each open set $U\supset C$ we have ${\sf tp}_{\sone(\opkfd,\op)}(X\setminus U) \le {\sf tp}_{\sone(\opkfd,\op)}(X)-2$. Then we have the contradiction that ${\sf tp}_{\sone(\opkfd,\op)}(X\setminus U) \le {\sf tp}_{\sone(\opkfd,\op)}(X)-1$.
$\square$

\begin{lemma}\label{finiteunion} For $X$ a metric space, and $Y_1,\cdots,Y_n$ subspaces such that for each $n$, ${\sf tp}_{\sone(\opkfd,\op)}(Y_n) \le\alpha$, also ${\sf tp}_{\sone(\opkfd,\op)}(\bigcup_{j\le n}Y_j) \le\alpha$.
\end{lemma}
{\flushleft{\bf Proof:}} For each $n$ we prove this by induction on $\alpha$. For $n=1$ there is nothing to prove. Thus assume $n>1$.
When $\alpha=0$ there is nothing to prove since each $Y_j$ then is compact and finite dimensional, as is their union. Thus, assume that $0<\alpha$ and we have proven this result for all $\beta<\alpha$. There are three cases to consider: $\alpha$ is finite, $\alpha$ is an infinite limit ordinal and $\alpha$ is an infinite successor ordinal.

First, the case when $\alpha$ is finite: Since we have already disposed of the case when $\alpha=0$, consider now the case where $\alpha = m+1$ and the result is proved for $\alpha\le m$. By Lemma \ref{twowinsinfinite} choose for $j\le n$ a compact finite dimensional set $C_j\subset Y_j$ such that for any open (in $Y_j$) set $U$ containing $C_j$, ${\sf tp}_{\sone(\opkfd,\op)}(Y_j\setminus U)\le m$.

Then $C = \cup_{j\le n}C_j$ is compact and finite dimensional, and for every open set $U\supset C$, for each $j$, ${\sf tp}_{\sone(\opkfd,\op)}(Y_j\setminus U)\le m$. Thus TWO can play as follows:
When ONE in inning 0 play an $O_0\in\opkfd$ for $\cup_{j\le n}Y_j$, TWO chooses $T_0\in O_0$ with $C\subseteq T_0$. Since now for $j\le n$ we have for $Z_j = Y_j\setminus T_0$ that ${\sf tp}_{\sone(\opkfd,\op)}(Z_j)\le m$, the induction hypothesis gives that ${\sf tp}_{\sone(\opkfd,\op)}(\cup_{j\le n}Z_j)\le m = {\sf tp}_{\sone(\opkfd,\op)}((\cup_{j\le n}Y_j)\setminus T_0)\le m$. It follows that ${\sf tp}_{\sone(\opkfd,\op)}(\cup_{j\le n}Y_j)\le m+1$, completing the induction step, and the proof for $\alpha$ finite.

In the case when $\alpha$ is a limit ordinal we can write $\alpha = \cup_{j\le m}A_j$ where each $A_j$ has order type $\alpha$, and in innings $\gamma\in A_j$, we use a winning strategy for TWO on $Y_j$ to choose an element $T_{\gamma}\in O_{\gamma}$. This plan produces a winning strategy in the $\alpha$-length game on $\cup_{j\le n}Y_j$.

In the case when $\alpha$ is an infinite successor ordinal, write $\alpha = \ell(\alpha)+n(\alpha)$ where $\ell(\alpha)$ is a limit ordinal and $n(\alpha)<\omega$. During the first $\ell(\alpha)$ innings we follow the plan above. After these innings for each $j$ the uncovered part $Z_j$ of $Y_j$ has  ${\sf tp}_{\sone(\opkfd,\op)}(Z_j)\le n(\alpha)$, and thus the uncovered part $Z$ of $\cup_{j\le n}Y_j$ has ${\sf tp}_{\sone(\opkfd,\op)}(Z)\le n(\alpha)$, by the finite case. This completes the proof.
$\square$\\

Theorem \ref{polspacerinfinity} below shows that Lemma \ref{finiteunion} fails for infinite unions.\\

{\flushleft 2.3 {\bf When ${\sf tp}_{\sone(\opkfd,\op)}(X)$ is finite and positive}}\\

Note that when $X$ is a compact finite-dimensional space, ${\sf tp}_{\sone(\opkfd,\op)}(X)=0$. Next we describe the structure of those metrizable $X$ with ${\sf tp}_{\sone(\opkfd,\op)}(X)$ finite.

Smirnov's compactum $\smirnovspace$ is constructed as follows: For positive integer $n$ define $S_n = \closedunitinterval^n$. Then define $\smirnovspace$ to be the one-point compactification of the topological sum $\oplus_{n=1}^{\infty} S_n$, say $\smirnovspace = (\oplus_{n=1}^{\infty} S_n)\bigcup\{p_{\omega}\}$. It is clear that $\smirnovspace$ is a union of countably many compact finite dimensional spaces.
\begin{lemma}\label{smirnovpowers} For each positive integer $n$, ${\sf tp}_{\sone(\opkfd,\op)}(\smirnovspace^n) = n$.
\end{lemma}
{\flushleft{\bf Proof:}} The proof is by induction on $n$.\\
\underline{$n=1$:}\\ 
Since $\smirnovspace$ is compact and not finite dimensional, $\opkfd\neq\emptyset$, and so TWO does not win in one inning. When ONE plays $O_0\in\opkfd$, then TWO chooses $T_0\in O_0$ such that $p_{\omega}\in T_0$. Then $\smirnovspace\setminus T_0$ is compact finite dimensional, and thus in the next inning TWO chooses $T_1\in O_1$ containing this set.\\
\underline{$n=k+1$:}\\
Assume that we have already proven the result for $n\le k$. Consider $n=k+1$: Player ONE will play open covers $O$ with the property that for any element $U$ of $O_1$ which contains the k+1 vector $(p_{\omega},\cdots,p_{\omega})$ there is some $m>k+1$ such that the projection of $U$ in each of the $k+1$ directions is disjoint from $\oplus_{j\le m}{\sf S}_j$. It follows that for these elements $U$ of ONE's move, $\smirnovspace^{k+1}\setminus U$ contains the closed subspace $\bigcup_{j=1}^{k+1}X_j$, where $X_j$ is the product $\prod_{i=1}^{k+1}Z_i$ where 
\[
  Z_i = \left\{\begin{tabular}{ll}
               $\smirnovspace$ & if $i\neq j$\\
               ${\sf S}_{k+1}$& otherwise 
               \end{tabular}
        \right.
\]
Each $X_j$ is homeomorphic to ${\sf S}_{k+1}\times \smirnovspace^{k}$. By Lemma \ref{cptfindimfactor}, ${\sf tp}_{\sone(\opkfd,\op)}({\sf S}_{k+1}\times \smirnovspace^k) = {\sf tp}_{\sone(\opkfd,\op)}(\smirnovspace^k) = k$. Then by Lemma \ref{finiteunion} ${\sf tp}_{\sone(\opkfd,\op)}(\smirnovspace^{k+1}\setminus U) \ge k$. Thus we conclude that ${\sf tp}_{\sone(\opkfd,\op)}(\smirnovspace^{k+1}) \ge k+1$. Proving that this inequality is in fact an equality is left to the reader. $\square$

{\flushleft 2.4 {{\bf When ${\sf tp}_{\sone(\opkfd,\op)}(X)$ is $\omega+1 $}}}\\

The following Lemma is useful for computing lowerbounds on ${\sf tp}_{\sone(\opkfd,\op)}(X)$:

\begin{lemma}\label{complementlemma} Let $X$ and $C$ be metrizable spaces with $X$ not countable dimensional. If $E\subset X\times C$ is countable dimensional, then there is a set $Y\subset X$ with $Y$ not countable dimensional and $E\bigcap (Y\times C) = \emptyset$.
\end{lemma} 
{\flushleft{\bf Proof:}} Suppose not. Then the set
\[
  Y:=\{x\in X: (\{x\}\times C) \bigcap E = \emptyset\}
\]
is countable dimensional, and so $Z = X\setminus Y$ is not countable dimensional. For each $x\in Z$ choose $\phi(x) \in C$ with $(x,\phi(x))\in E$. Then the set $T:=\{(x,\phi(x)):x\in Z\}$, being a subset of $E$, is countable dimensional. But the function $f:T\longrightarrow Z$ defined by $f(x,y) = x$, the projection on the first coordinate, being the restriction of $\Pi_1$ to $T$, is 
continuous, and one-to-one and onto $Z$. Thus, the inverse of $f$ is a both open and closed function from $Z$ to $T$. Since $T$ is countable dimensional, a theorem of Arhangel'ski\v{i} (see Theorem 7.7 in \cite{engelkingpol}) 
implies that $Z$ is countable dimensional, a contradiction. It follows that $Y$ is not countable dimensional.
$\square$

Pol's compactum $\polspace$ is constructed as follows: Start with a complete hereditarily strongly infinite dimensional totally disconnected metric space $M$ and then compactify it such that the extension $L$ is a union of countably many compact, finite dimensional spaces. Then $\polspace = M\bigcup L$. For the rest of the paper fix a representation of $L$ as a union of countably many compact finite dimensional sets, say $L = \bigcup_{n=1}^{\infty}L_n$, where for $m<n$ we have $L_m\subset L_n$ and $dim(L_m)<dim(L_n)$.

\begin{theorem}\label{polspace} For each positive integer $n$, ${\sf tp}_{\sone(\opkfd,\op)}(\polspace^n) = \omega\cdot n +1$
\end{theorem}
{\flushleft{\bf Proof:}} We use induction on $n$.\\
\underline{n=1:}\\ 
\underline{${\sf tp}_{\sone(\opkfd,\op)}(\polspace) \le \omega+1$:} In inning $n$ when ONE plays the open cover $O_n\in\opkfd$, TWO chooses $T_n\in O_n$ with $L_n\subset T_n$. After $\omega$ innings $L\subset \bigcup_{n=1}^{\infty} T_n$ and the part of $\polspace$ not yet covered by TWO is a compact set contained in the totally disconnected space $M$, and thus is a compact zerodimensional set. In inning $\omega+1$ TWO chooses $T_{\omega}\in O_{\omega}$ containing this compact zerodimensional set. \\
\underline{${\sf tp}_{\sone(\opkfd,\op)}(\polspace) \ge \omega+1$:} $\polspace$ is not countable dimensional. By Lemma \ref{opomegalength}, ${\sf tp}_{\sone(\opkfd,\op)}(\polspace)>\omega$.
{\flushleft{\underline{$n>1$:}}} Assume that the statement is true for $k<n$.\\ 
\underline{${\sf tp}_{\sone(\opkfd,\op)}(\polspace^n)\le \omega\cdot n + 1$:} For $j<n$ and $m<\omega$ define $\polspace(j,m) = \polspace^j \times L_m\times \polspace^{n-j-1}$. For each such $(j,m)$: $\polspace(j,m)$ is homeomorphic to $L_m\times\polspace^{n-1}$ and by Lemma \ref{cptfindimfactor} and the induction hypothesis ${\sf tp}_{\sone(\opkfd,\op)}(\polspace(j,m)) = \omega\cdot(n-1)+1$. Let $F_{j,m}$ be TWO's winning strategy in $\gone^{\omega\cdot(n-1) + 1}(\opkfd,\op)$ on $\polspace(j,m)$.

Now write $\omega\cdot(n-1)$ as a union $\bigcup_{m<\omega}\bigcup_{j<n}S_{j,m}$ where each $S_{j,m}$ has order type $\omega\cdot(n-1)$ and for $(i,m)\neq (j,\ell)$, $S_{i,m}\cap S_{j,\ell}=\emptyset$. Enumerate each $S_{j,m}$ in an order preserving way as $(s^{j,m}_{\gamma}:\gamma<\omega\cdot(n-1))$. 

In the first $\omega\cdot(n-1)$ innings, in inning $\gamma$, when ONE plays open cover $O_{\gamma}$ of $\polspace^n$, TWO responds as follows: Choose ${j,m}$ with $\gamma\in S_{j,m}$, and choose $\delta$ with $\gamma = s^{j,m}_{\delta}$. Write $O^{j,m}_{\delta}$ for $O_{\gamma}$. Then TWO applies the winning strategy for inning $\delta$ of the game on $\polspace(j,m)$ to choose $T_{\gamma}\in O_{\gamma}$ by $T_{\gamma} = F_{j,m}(O^{j,m}_{\nu}:\nu\le \delta)$. Then after $\omega\cdot(n-1)$ innings TWO has constructed the open set $V=\bigcup_{\gamma<\omega\cdot(n-1)}T_{\gamma}$ in such a way that for each $(j,m)$ the set $C(j,m) = \polspace(j,m)\setminus V$ is compact and finite dimensional.

Next enumerate $\{(j,m):j<n,\, m<\omega\}$ bijectively as $\{(j_k,m_k):k<\omega\}$. In inning $\omega\cdot(n-1)+k$ TWO chooses $T_{\omega\cdot(n-1)+k}\in O_{\omega\cdot(n-1)+k}$ containing $C(j_k,m_k)$. Put $U = \cup_{\gamma<\omega\cdot n}T_{\gamma}$. 

Note that $M^n\setminus U = \polspace^n\setminus U$: For if $(x_1,\cdots,x_n)\not\in M^n$, choose a $j$ with $x_j\in L$, and choose an $m$ with $x_j\in L_m$. Then $(x_1,\cdots,x_n)\in \polspace(j,m)$ and thus if it is not in $V$, it is in $C(j,m)$ and thus in $U$. Since $M^n\setminus U$ is compact and totally disconnected it is  compact and zero-dimensional, and so in one more inning TWO covers $M^n\setminus U$.

{\flushleft{\underline{${\sf tp}_{\sone(\opkfd,\op)}(\polspace^n)\ge \omega\cdot n + 1$:}}} Let $\alpha = {\sf tp}_{\sone(\opkfd,\op)}(\polspace^n)$. Let $F$ be a winning strategy for TWO in the game $\gone^{\alpha}(\opkfd,\op)$. Put $T = L^n$, a compactly countable dimensional set. By Theorem \ref{twoskerneldirectedcofinal} choose a set $S\subset \polspace^n$ which is a countable union of compact finite dimensional sets such that $T\subseteq S$ and for any closed set $C\subset (\polspace^n\setminus S)$ there is an $\omega$-length $F$-play $O_0,\, T_0,\, \cdots,\, O_n,\, T_n\,\cdots$ such that $T\subseteq \bigcup_{n=1}^{\infty}T_n \subseteq X\setminus C$. Since $\polspace$ is not countable dimensional, choose by Lemma \ref{complementlemma} a point $x\in \polspace$ such that $\{x\}\times\polspace^{n-1} \cap S = \emptyset$. Now $C=\{x\}\times\polspace^{n-1}$ is closed in $\polspace^n$, and disjoint from $S$. Choose an $\omega$-length play  $O_0,\, T_0,\, \cdots,\, O_n,\, T_n\,\cdots$ where TWO used $F$ and $\{x\}\times\polspace^{n-1} \cap (\bigcup_{n<\infty} T_n)=\emptyset$. Since $\{x\}\times\polspace^{n-1}$ is homeomorphic to $\polspace^{n-1}$ and by the induction hyptohesis ${\sf tp}_{\sone(\opkfd,\op)}(\polspace^{n-1})=\omega\cdot(n-1)+1$, it requires at least $\omega\cdot(n-1)+1$ more innings for TWO to win. 
$\square$

In particular for each $n$, $\polspace^n$ also satisfies $\sone(\opkfd,\op)$. 
\begin{corollary}\label{polsums}
For each positive integer $n$, ${\sf tp}_{\sone(\opkfd,\op)}(\polspace^n) = \sum_{j=1}^n {\sf tp}_{\sone(\opkfd,\op)}(\polspace)$.
\end{corollary}

\begin{lemma}\label{polspacerinfinity} Let $X$ be metrizable with ${\sf tp}_{\sone(\opkfd,\op)}(X) = \omega$. For each $m$, 
\[
  {\sf tp}_{\sone(\opkfd,\op)}(\polspace^m\times X) = {\sf tp}_{\sone(\opkfd,\op)}(\polspace^m) + {\sf tp}_{\sone(\opkfd,\op)}(X)=\omega\cdot (m+1).
\]
\end{lemma}
{\flushleft {\bf Proof:}} For each $n$ Let $X_n$ be compact finite dimensional so that $X=\cup_{n<\omega}X_n$. We prove this theorem by induction on $m$.
{\flushleft{The case when $m=1$:}}
{\flushleft{\underline{(1) ${\sf tp}_{\sone(\opkfd,\op)}(\polspace\times X)\le \omega\cdot 2$:}}}\\
By Lemma \ref{cptfindimfactor} and Theorem \ref{polspace}, for each $n$ TWO has a winning strategy $F_n$ in $\gone^{\omega+1}(\opkfd,\op)$ on $\polspace\times X_n$. 

Write $\omega=\cup_{n<\infty}S_n$ where each $S_n$ is infinite, and $S_m\cap S_n=\emptyset$ for $m\neq n$. For each $n$ we enumerate $S_n$ in order as $(s^n_j:j<\omega)$. 

In the first $\omega$ innings, when ONE plays $O_k$ in inning $k$, TWO chooses $n$ with $k\in S_n$ and then fixes $j$ with $k=s^n_j$, and think of $O_k$ as $O^n_j$, the $j$-th move of ONE in the game on $\polspace\times X_n$. Then TWO chooses $T_k\in O_k$ using $F_n(O^n_0,\cdots,O^n_j)$. 

Let $U$ be the union of TWO's moves during the first $\omega$ innings. For each $n$, since TWO was using the winning strategy $F_n$ through $\omega$ innings on $\polspace\times X_n$, the set $C_n=(\polspace\times X_n)\setminus U$ is compact finite dimensional. TWO covers these countably many compact finite dimensional sets in the next $\omega$ innings. 

{\flushleft{\underline{(2) ${\sf tp}_{\sone(\opkfd,\op)}(\polspace\times X)\ge \omega\cdot 2$:}}}\\
Put $\alpha={\sf tp}_{\sone(\opkfd,\op)}(\polspace\times X)$. Let $F$ be a winning strategy for TWO in $\gone^{\alpha}(\opkfd,\op)$ played on $\polspace\times X$. The set $T = L\times X$ is compactly countable dimensional. By Theorem \ref{twoskerneldirected} choose a countable dimensional set $S\subset\polspace\times X$  
such that $T\subseteq S$ and for any closed set $C\subset (\polspace\times X)\setminus S$ there is an $\omega$-length $F$-play $O_0,\, T_0,\, \cdots,\, O_n,\, T_n,\, \cdots$ such that $T\subset \bigcup_{n=1}^{\infty}T_n\subset (\polspace\times X)\setminus C$. By Lemma \ref{complementlemma} choose an $x\in M$ with the closed set $C = \{x\}\times X$ disjoint from $S$. It follows that in some $F$-play, after $\omega$ innings TWO still has not covered any point in $C$, and thus the game will last at least another $\omega$ innings. Thus $\alpha\ge\omega\cdot 2$.

This completes the proof when $m=1$. Now assume that $m>1$ and the theorem holds for all $j<m$:
{\flushleft{\underline{(1) ${\sf tp}_{\sone(\opkfd,\op)}(\polspace^m\times X)\le \omega\cdot (m+1)$:}}}\\ 
By Lemma \ref{cptfindimfactor} and Theorem \ref{polspace}, for each $n$ TWO has a winning strategy $F_n$ in $\gone^{\omega\cdot m + 1}(\opkfd,\op)$ on $\polspace^{m}\times X_n$. 

Write $\omega\cdot m=\cup_{n<\infty}S_n$ where each $S_n$ is infinite, and $S_m\cap S_n=\emptyset$ for $m\neq n$. For each $n$ we enumerate $S_n$ in order as $(s^n_{\gamma}:\gamma<\omega\cdot m)$. 

In the first $\omega\cdot m$ innings, when ONE plays $O_{\gamma}$ in inning $\gamma$, TWO chooses $n$ with $\gamma\in S_n$ and then fixes $\delta$ with $\gamma=s^n_{\delta}$, and think of $O_{\gamma}$ as $O^n_{\delta}$, the $\delta$-th move of ONE in the game on $\polspace^m\times X_n$. Then TWO chooses $T_{\gamma}\in O_{\gamma}$ using $F_n(O^n_0,\cdots,O^n_{\delta})$.

Let $U$ be the union of TWO's moves during the first $\omega\cdot m$ innings. For each $n$, since TWO was using the winning strategy $F_n$ through $\omega\cdot m$ innings on $\polspace^m\times X_n$, the set $C_n=(\polspace^m\times X_n)\setminus U$ is compact finite dimensional. TWO covers these countably many compact finite dimensional sets in the next $\omega$ innings. 

{\flushleft{\underline{(2) ${\sf tp}_{\sone(\opkfd,\op)}(\polspace^m\times X)\ge \omega\cdot (m+1)$:}}}\\ 
Let $F$ be a winning strategy for TWO in $\gone^{\alpha}(\opkfd,\op)$ played on $\polspace^m\times X$. The set $T = L^m\times X$ is compactly countable dimensional. By Theorem \ref{twoskerneldirected} choose a countable dimensional set $S\subset \polspace^m\times X$ such that $T\subseteq S$ and for each closed set $C\subset (\polspace^m\times X)\setminus S$ there is an $\omega$-length $F$-play $O_0,\, T_0,\, \cdots,\, O_n,\, T_n,\, \cdots$ with $T\subset\cup_{n<\infty}T_n\subset (\polspace^m\times X)\setminus C$. Then as $\polspace$ is not countable dimensional, Lemma \ref{complementlemma} implies that there is an $x\in\polspace$ with the closed set $C = \{x\}\times (\polspace^{m-1}\times X)$ disjoint from $S$. Then choose an $\omega$-length play $O_0,\, T_0,\, \cdots,\, O_n,\, T_n,\, \cdots$ with $T\subset\cup_{n<\infty}T_n\subset (\polspace^m\times X)\setminus C$. Since $C = \{x\}\times (\polspace^{m-1}\times X)$ is disjoint from the set covered by TWO in these innings, and is homeomorphic to $\polspace^{m-1}\times X$ the induction hypothesis implies that at least $\omega\cdot m$ additional inning are needed for TWO to cover $C$. It follows that $\alpha\ge \omega + \omega\cdot m = \omega\cdot (m+1)$. 

This completes the induction step and the proof.
$\square$\\

\begin{theorem}\label{successorproductwithKCD} Let $X$ be a metrizable space with ${\sf tp}_{\sone(\opkfd,\op)}(X)$ a successor ordinal. Let $Y$ be a metric space with ${\sf tp}_{\sone(\opkfd,\op)}(Y)\le \omega$. Then:
\[
  {\sf tp}_{\sone(\opkfd,\op)}(X\times Y) \le {\sf tp}_{\sone(\opkfd,\op)}(X) + \omega.
\]
\end{theorem}
{\flushleft {\bf Proof:}} 

Write $\alpha = \beta+m$ where $\beta<\alpha={\sf tp}_{\sone(\opkfd,\op)}(X)$ is a limit ordinal, and $0<m<\omega$. Write $\beta=\cup_{n<\infty}S_n$ where each $S_n$ is a subset of $\beta$ of order type $\beta$, and $S_m\cap S_n=\emptyset$ for $m\neq n$. For each $n$ we enumerate $S_n$ as $(s^n_{\gamma}:\gamma<\beta)$. 

When ONE plays $O_{\gamma}$ in inning $\gamma<\beta$, TWO chooses $n$ with $\gamma\in S_n$ and then fixes $\nu$ with $\gamma=s^n_{\nu}$, and think of $O_{\gamma}$ as $O^n_{\nu}$, the $\nu$-th move of ONE in the game on $X\times Y_n$. Then TWO chooses $T_{\gamma}\in O_{\gamma}$ using $F_n(O^n_0,\cdots,O^n_{\nu})$. 

Let $U$ be the union of the moves TWO made during these $\beta$ innings. For each $n$, since TWO was using the winning strategy $F_n$ through $\beta$ innings on $X\times Y_n$, we have ${\sf tp}_{\sone(\opkfd,\op)}((X\times Y_n)\setminus U)\le m$. But then, by Proposition \ref{omegalength}, for each $n$ the set $C_n = (X\times Y_n)\setminus U$ is a union of countably many compact finite dimensional sets. During the next $\omega$ innings TWO covers these. Thus we have that ${\sf tp}_{\sone(\opkfd,\op)}(X\times Y)\le \alpha+\omega$. $\square$

Theorem \ref{polspacerinfinity} below shows that the value of ${\sf tp}_{\sone(\opkfd,\op)}(X\times Y)$ obtained in Theorem \ref{successorproductwithKCD} cannot be improved. We suspect that the upper bound in Theorem \ref{successorproductwithKCD} is in fact achieved:
\begin{problem}
If $X$ be a metrizable space with ${\sf tp}_{\sone(\opkfd,\op)}(X)$ a successor ordinal and $Y$ is a metric space with ${\sf tp}_{\sone(\opkfd,\op)}(Y)\le \omega$, does it follow that:

\[
  {\sf tp}_{\sone(\opkfd,\op)}(X\times Y) = {\sf tp}_{\sone(\opkfd,\op)}(X) + \omega?
\]
\end{problem}

{\flushleft 2.4 {{\bf When ${\sf tp}_{\sone(\opkfd,\op)}(X)$ is a limit ordinal}}}\\

\begin{theorem}\label{productwithKCD} Let $X$ be a metrizable space with ${\sf tp}_{\sone(\opkfd,\op)}(X)=\alpha$ a limit ordinal, and let $Y$ be a metric space with ${\sf tp}_{\sone(\opkfd,\op)}(Y)\le \omega$. Then:
${\sf tp}_{\sone(\opkfd,\op)}(X\times Y) = \alpha$.
\end{theorem}
{\flushleft {\bf Proof:}} Since ${\sf tp}_{\sone(\opkfd,\op)}(Y)\le \omega$, by Proposition \ref{omegalength} choose compact finite dimensional sets $Y_n\subset Y$ with $Y = \cup_{n<\omega}Y_n$. By Lemma \ref{cptfindimfactor} fix for each $n$ a winning strategy $F_n$ for TWO in the game $\gone^{\alpha}(\opkfd,\op)$ played on $X\times Y_n$.

By Lemma \ref{opvsom} (3), since $X$ is homeomorphic to a closed subspace of $X\times Y$, we have $\alpha\le {\sf tp}_{\sone(\opkfd,\op)}(X\times Y)$. We must show that ${\sf tp}_{\sone(\opkfd,\op)}(X\times Y)\le \alpha$.

Write $\alpha=\cup_{n<\infty}S_n$ where each $S_n$ is a subset of $\alpha$ of order type $\alpha$, and $S_m\cap S_n=\emptyset$ for $m\neq n$. For each $n$ we enumerate $S_n$ in order as $(s^n_{\gamma}:\gamma<\alpha)$. 

When ONE plays $O_{\gamma}$ in inning $\gamma<\alpha$, TWO chooses $n$ with $\gamma\in S_n$ and then fixes $\nu$ with $\gamma=s^n_{\nu}$, and think of $O_{\gamma}$ as $O^n_{\nu}$, the $\nu$-th move of ONE in the game on $X\times Y_n$. Then TWO chooses $T_{\gamma}\in O_{\gamma}$ using $F_n(O^n_0,\cdots,O^n_{\nu})$. 

Let $U$ be the union of the moves TWO made during these $\alpha$ innings. For each $n$, since TWO was using the winning strategy $F_n$ through $\alpha$ innings on $X\times Y_n$, we have $X\times Y_n\subseteq U$. Thus, TWO won in $\alpha$ innings. This completes the proof. $\square$\\

Let $X$ be $\oplus_{n=1}^{\infty}\polspace^n$, the topological sum of the finite powers of the Pol compactum. $\omega^2= {\sf tp}_{\sone(\opkfd,\op)}(X)$. This is a locally compact space. Let $\polspace_\omega = (\oplus_{n=1}^{\infty}\polspace^n)\cup\{p_{\omega}\}$ be the one-point compactification of $X$. It is easy to see that ${\sf tp}_{\sone(\opkfd,\op)}(\polspace_{\omega})=\omega^2$. \\

We now show that ${\sf tp}_{\sone(\opkfd,\op)}(X\times\polspace)=\omega^2$.

\begin{lemma}\label{polomegaproducts} ${\sf tp}_{\sone(\opkfd,\op)}(\polspace_{\omega}\times \polspace) = \omega^2$. Thus, ${\sf tp}_{\sone(\opkfd,\op)}(\polspace_{\omega}\times X) = \omega^2$.
\end{lemma}
{\flushleft{\bf Proof:}}  During the first $\omega+1$ innings, cover $\{p_{\omega}\}\times\polspace$. The remaining part of the space is a closed subset of $X$ and thus requires at most $\omega^2$ more innings. Thus $\omega+1 + \omega^2 \ge {\sf tp}_{\sone(\opkfd,\op)}(\polspace_{\omega}\times\polspace)$. But $\polspace_{\omega}$ is a closed subset of $\polspace_{\omega}\times\polspace$, and so also $\omega^2 \le {\sf tp}_{\sone(\opkfd,\op)}(\polspace_{\omega}\times\polspace)$. $\square$

\section{The ordinal ${\sf tp}_{\Sc(\op,\op)}(S)$ and the game dimension.}

We now define an ordinal-valued function on the set of subspaces of the Hilbert cube such that this function
\begin{itemize}
\item{coincides with Lebesgue covering dimension in the case of subspaces with finite covering dimension, and}
\item{captures an important aspect of extending the covering dimension to subspaces that are not finite dimensional.}
\end{itemize} 
Other important examples of such ordinal-valued functions were developed by P. Borst \cite{borstcspaces}, and R. Pol \cite{rpol}. 

Let $\alpha$ be an ordinal number. The game $\Gc^{\alpha}(\mathcal{A},\mathcal{B})$ associated with the selection principle $\Sc(\mathcal{A},\mathcal{B})$ is as follows: The players play an inning per $\gamma<\alpha$. In the $\gamma$-th inning ONE first chooses an $A_{\gamma}\in\mathcal{A}$: TWO then responds with a  pairwise disjoint family of sets $B_{\gamma}$ that refines $A_{\gamma}$. A play $A_0,\, B_0,\, \cdots,\, A_{\gamma},\, B_{\gamma},\, \cdots$ of length $\alpha$ is won by TWO if $\bigcup_{n\in\naturals}B_n\in\mathcal{B}$. Else, ONE wins.

When for a set S and families $\mathcal{A}$ and $\mathcal{B}$ there is an ordinal number $\alpha$ such that TWO has a winning strategy in the game $\Gc^{\alpha}(\mathcal{A},\mathcal{B})$ played on S, then we define:
\[
  {\sf tp}_{\Sc(\mathcal{A},\mathcal{B})}(S) = \min\{\alpha:\mbox{ TWO has a winning strategy in }\Gc^{\alpha}(\mathcal{A},\mathcal{B}) \mbox{ on }S\}.
\]

For  $X$ a separable finite-dimensional metric space, $\dim(X)$ denotes the Lebesgue covering dimension of $X$. 
The starting point for this exploration is the following game-theoretic result:
\begin{lemma}[\cite{lbscgame}]\label{findimgc} Let $X$ be a separable metric space and let n be a nonnegative integer. The following are equivalent:
\begin{enumerate}
\item{$\dim$(X) = n.}
\item{${\sf tp}_{\Sc(\op,\op)}(X)=n+1$.}
\end{enumerate}
\end{lemma}

We define for separable metric space $X$ the \emph{game dimension} of $X$, denoted $\dim_G(X)$, by 
\[
  1+ \dim_G(X) = {\sf tp}_{\Sc(\op,\op)}(X).
\]
The for $X$ a finite-dimensional separable metric space, $\dim_G(X)=\dim(X)$. By a theorem of Nagami \cite{nagami} and Smirnov  \cite{smirnov} (independently), each separable metric space is the union of $\le\omega_1$ zero-dimensional subsets. This implies
\begin{theorem}\label{existence}
For each separable metric space $X$, $\dim_G(X) \le \omega_1$.
\end{theorem}
The following three results are easy to prove.
\begin{theorem}\label{cspaces}
If $X$ is a separable metric space for which $\dim_G(X)<\omega_1$, then $X$ has property $\Sc(\op,\op)$.
\end{theorem}
Since the Hilbert cube $\hilbertcube = \lbrack 0,\, 1\rbrack^{\naturals}$ does not have the property $\Sc(\op,\op)$, we have $dim_G(\hilbertcube)=\omega_1$.

\begin{lemma}[Subspace Lemma]\label{subspace}
Let $Y$ be a closed subspace of $X$. Then $\dim_G(Y) \le \dim_G(X)$.
\end{lemma}

\begin{lemma}[Addition Lemma]\label{addition}
Let $Y\subseteq X$ be such that $\dim_G(Y) = \beta$, and $\alpha$ is minimal such that for each closed set $C\subset X\setminus Y$ we have $\dim_G(C) \le \alpha$. Then $\dim_G(X) = \beta+\alpha$.
\end{lemma}

We have the following particularly satisfying property of game dimension in the case of countable dimensional spaces:
\begin{lemma}[\cite{lbscgame}]\label{ctbldimgc} Let $X$ be a separable metric space. The following are equivalent:
\begin{enumerate}
\item{X is countable dimensional.}
\item{$\dim_G(X)=\omega$.}
\end{enumerate}
\end{lemma}

Now consider the behavior of game dimension on spaces $X$ with $\omega< \dim_G(X) <\omega_1$. As noted in Theorem \ref{cspaces} these spaces are among the spaces with property $\Sc(\op,\op)$. A number of important subclasses of the spaces with property $\Sc(\op,\op)$ have been identified and play an important role in developing game dimension. The following results point out some of these connections. 

\begin{theorem}\label{scequivs} $\Sc(\op,\op)  = \Sc(\opkfd,\op) =\Sc(\opfd,\op)$.
\end{theorem}
{\bf Proof:} We show that $ \Sc(\opfd,\op) \Rightarrow \Sc(\op,\op)$. Let $X$ satisfy $\Sc(\opfd,\op)$. Let $(\mathcal{U}_n:n<\infty)$ be a sequence of open covers of $X$.  Now we use the technique in \cite{lbhaver} to convert this sequence of open covers to a sequence of {\sf FD}-covers: Write $\naturals = \bigcup_{k<\infty} Y_k$ where each $Y_k$ is infinite and the $Y_k$'s are pairwise disjoint.

Fix a $k$, a finite dimensional subset $C$ of the space, and a set $I\subset Y_k$ with $|I| = dim(C)+1$. Then for each $i\in I$ $\mathcal{U}_i$ is an open cover of $C$ and we can find a pairwise disjoint refinement $\mathcal{V}_i$ of $\mathcal{U}_i$, where $\mathcal{V}_i$ consists of sets open in the relative topology of $C$, such that $\bigcup_{i\in I}\mathcal{V}_i$ is an open cover of $C$. Without loss of generality we may assume that the sets in $\mathcal{V}_i$ are open in $X$, and pairwise disjoint in $X$. Define
\[
  U(C,I) = \bigcup(\bigcup_{i\in I}\mathcal{V}_i).
\]
Finally, set $\mathcal{V}_k = \{U(C,I): C\subset X \mbox{ finite dimensional and } I\subset Y_k \mbox{ with } |I| = dim(C)+1\}$. Then each $\mathcal{V}_k$ is in $\opfd$. 

Applying $\Sc(\opfd,\op)$ to the sequence $(\mathcal{V}_k:k<\infty)$ we find a sequence $(\mathcal{H}_k:k<\infty)$ so that each $\mathcal{H}_k$ is pairwise disjoint, refines $\mathcal{V}_k$, and $\bigcup_{k<\infty}\mathcal{H}_k$ is an open cover of $X$. For each $k$, and for each $W\in\mathcal{H}_k$, choose a finite dimensional set $C_W$ and a finite set $I_W\subset Y_k$ so that $W\subset U(C_W,I_W)$. For each $i\in I_W$, choose $\mathcal{V}_i(W)$ an open pairwise disjoint refinement of $\mathcal{U}_i$ so that $U(C_W,I_W) = \bigcup(\bigcup_{i\in I_W}\mathcal{V}_i(W))$. Put $I_k = \bigcup_{W\in\mathcal{H}_k}I_W$, a subset of $Y_k$. For each $i\in I_k$ define 
\[
  \mathcal{G}_i = \{W\cap V: W\in\mathcal{H}_i \mbox{ and }V\in \mathcal{V}_i(W)\}.
\]
Each $\mathcal{G}_i$ is a pairwise disjoint open refinement of the corresponding $\mathcal{U}_i$, and the union of all the $\mathcal{H}_i$'s is an open cover of $X$. 
$\square$

The inequalities in the following theorem can be proved also for the family $\opkfd$ substituting for $\opfd$. However, the result given here seems optimal.

\begin{theorem}\label{scvsopfd} Let $X$ be an infinite dimensional metrizable space. Then: 
   \[
     \dim_G(X)\le \left\{
        \begin{tabular}{ll}
                      ${\sf tp}_{\sone(\opfd,\op)}(X)$        & \mbox{ if ${\sf tp}_{\sone(\opfd,\op)}(X)$ is a limit ordinal}\\
                      ${\sf tp}_{\sone(\opfd,\op)}(X)+\omega$ & \mbox{ otherwise.}\\
        \end{tabular}
\right.
\]
\end{theorem}
{\flushleft{\bf Proof:}} Let $\alpha$ be an ordinal such that TWO has a winning strategy in the game $\gone^{\alpha}(\opfd,\op)$. Let $F$ be a winning strategy for TWO. We shall use $F$ to define a winning strategy for TWO in $\Gc^{\beta}(\op,\op)$ for the appropriate $\beta$. We shall repeatedly use Lemma \ref{setandopenset} and Lemma \ref{findimgc}. We may assume without loss of generality that the open covers played by player ONE of $\Gc^{\beta}(\op,\op)$ never include $X$ as an element.

To begin, choose by Lemma \ref{setandopenset} a finite dimensional set $C_0\subset X$ such that there is for each proper open set $U\supseteq C_0$ an $O\in\opfd$ with $F(O)=U$. Put $n_0 = dim(C_0)+1$. During the first $n_0$ innings of $\Gc^{\beta}(\op,\op)$, TWO chooses for each $O_j\in\op$ a disjoint open refinement $T_j$ such that $C_0\subseteq W_0 = \bigcup(\bigcup_{j\le n_0}T_j)$ (this uses Lemma \ref{findimgc}). Since $X$ is infinite dimensional $X\neq W_0$. Then choose a $B_0\in\opfd$ with $W_0=F(B_0)$. 

For the second stage of the game, choose by Lemma \ref{setandopenset} a finite dimensional set $C_1\subset X$ such that there is for each proper open set $U\supseteq C_1$ a $\mathcal{U}\in\opfd$ with $U=F(B_1,\mathcal{U})$. Put $n_1 = dim(C_1)+1$. During the next $n_1$ innings of $\Gc^{\beta}(\op,\op)$, TWO chooses for each $O_j\in\op$ a disjoint open refinement $T_j$ of $O_j$ such that $C_1\subseteq W_1 = \bigcup(\bigcup_{j=n_0+1}^{n_0+n_1}T_j)$ (this uses Lemma \ref{findimgc}). As before $X\neq W_1$: Choose $B_1\in\opfd$ with $F(B_0,B_1)=W_1$.

Continuing like this we find that after $\omega$ innings we have chosen finite dimensional subsets $C_0,\cdots,C_j,\cdots$ of $X$, open disjoint refinements $T_j$ of open covers $O_j$, and nonnegative integers $n_0,\cdots, n_j,\cdots$, $j<\omega$, and elements $B_j$ of $\opfd$ such that for all $j\ge 0$:
\begin{enumerate}
\item{$n_j = dim(C_j)+1$;}
\item{$C_j\subseteq \bigcup(\bigcup_{i=n_{j-1}+1)}^{n_{j-1}+n_j} T_i) = W_j \neq X$;}
\item{$F(B_0,\cdots,B_j) = W_j$.}
\end{enumerate}

Assume that $\nu<\alpha$ is an infinite ordinal and that for each $\rho<\nu$ we have selected a finite dimensional set $C_{\rho}$, an element $B_{\rho}$ of $\opfd$, an element $W_{\rho}$ of $B_{\rho}$, an open cover $O_{\rho}$ and a disjoint open refinement $T_{\rho}$ of $O_{\rho}$ such that the following hold: Write $\rho$ as 
\[
  \rho = \omega^{\beta_1}\cdot n_1 + \cdots + \omega^{\beta_k}\cdot n_k
\]
where now $\beta_1> \cdots >\beta_k\ge 0$, and each $n_i$ is positive.
\begin{enumerate}
\item{$C_{\rho}$ is such that for each proper open set $U\supseteq C_{\rho}$ there is a $B\in\opfd$ with $U = F((B_{\sigma}:\sigma<\rho)\frown(B))$;}
\item{If $\beta_k=0$, put $\mu=\sup\{\sigma<\rho:\, \sigma\mbox{ a limit ordinal}\}$, so that $\rho = \mu+n_k$. And set $m = \sum_{\mu\le j<\mu+n_k}(dim(C_{j})+1)$. Then we already have available $O_{\mu},\, \cdots,\, O_{\mu+m+dim(C_{\rho})+1}$ and $T_{\mu},\cdots, T_{\mu+m+dim(C_{\rho})+1}$, and $S_{\rho}\subseteq \bigcup(\bigcup_{i=\mu+m+1}^{\mu+m+dim(C_{\rho})+1} T_i)=W_{\rho}$;}
\item{If $\beta_k>0$, then $\rho$ is a limit ordinal. And set $m = dim(C_{\rho})+1$. Then we already have available $O_{\rho+1},\, \cdots,\, O_{\rho+m}$ and $T_{\rho+1},\cdots, T_{\rho+m}$, and $S_{\rho}\subseteq \bigcup(\bigcup_{i=\rho+1}^{\rho+m} T_i)=W_{\rho}$.}
\item{$W_{\rho}=F(B_{\sigma}:\sigma\le\rho)$.}
\end{enumerate}
We now describe what happens at $\nu$. \\
{\flushleft{\bf Case 1:} $\nu$ is a limit ordinal:}
By Lemma \ref{setandopenset} choose a finite dimensional subset $C_{\nu}$ of $X$ such that for each proper open set $U\supseteq C_{\nu}$ there is a $\mathcal{U}\in\opfd$ with $U = F((B_{\rho}:\rho<\nu)\frown(\mathcal{U}))$ and put $m = dim(C_{\nu})$. In innings $j\in\{\nu,\, \nu+1,\, \cdots, \nu+m\}$ of $\Gc^{\beta}(\op,\op)$, choose open disjoint refinements $T_j$ of $O_j$ such that $C_{\nu}\subseteq \bigcup(\bigcup_{j=\nu}^{\nu+m}T_j) = W_{\nu}$ (this uses Lemma \ref{findimgc}). Since $X$ is infinite dimensional, $W_{\nu}\neq X$, so choose by Lemma \ref{setandopenset} a $B_{\nu}\in\opfd$ with $W_{\nu} = F(B_{\rho}:\rho\le \nu)$.\\
{\flushleft{\bf Case 2:} $\nu$ is a successor ordinal:}
Write $\nu = \rho + k+1$ for some nonnegative integer $k$ and limit ordinal $\rho$. By the induction hypothesis we have finite dimensional sets $C_{\rho}, \cdots, C_{\rho+k}$ available, and with $\ell = \sum_{j=\rho}^{\rho+k}(dim(C_j)+1)$, we also have available the disjoint open refinements $T_{\rho},\cdots, T_{\rho+\ell}$ of the open covers $O_{\rho},\cdots,O_{\rho+\ell}$ of $X$. Now choose by Lemma \ref{setandopenset} a finite dimensional subset $C_{\nu}$ of $X$ such that for each proper open set $U\supseteq C_{\nu}$ there is a $\mathcal{U}\in\opfd$ such that $U = F((B_{\sigma}:\sigma<\nu)\frown(\mathcal{U}))$. Put $m=dim(C_{\nu})+1$. During the next $m$ innings of $\Gc^{\beta}(\op,\op)$ TWO chooses open disjoint refinements $T_{\rho+\ell+i}$ of $O_{\rho+\ell+i}$, $1\le i\le m$, so that $C_{\nu}\subseteq \bigcup(\bigcup_{i=\rho+\ell+1}^{\rho+\ell+m}T_i) = W_{\nu} \neq X$ (this uses Lemma \ref{findimgc}), and then choose $B_{\nu}\in\opfd$ such that $W_{\nu} = F(B_{\sigma}:\sigma\le \nu)$.

Thus, we see that the recursive conditions hold also at $\nu$. Thus for a play $((O_{\sigma}, T_{\sigma}):\sigma<\beta)$ of $\Gc^{\beta}(\op,\op)$ played according to this strategy we find a sequence $((B_{\nu},W_{\nu}):\nu<\alpha)$ which is an $F$-play of $\gone^{\alpha}(\opfd,\op)$, thus won by TWO, for which the union of the set of $W_{\nu}$ is covered by the union of the set of $T_{\sigma}$'s. Moreover, the strategy described shows that if $\alpha$ is a limit ordinal, then $\beta=\alpha$ works, and if $\alpha$ is a successor ordinal, then TWO wins $\Gc^{\beta}(\op,\op)$ at some inning before $\alpha+\omega$.
$\square$

Note that since ${\sf tp}_{\sone(\opfd,\op)}(X)\le {\sf tp}_{\sone(\opkfd,\op)}(X)$ we find that 
\begin{corollary}\label{scvsopkfd} For $X$ an infinite dimensional separable metric space,
\[
  \dim_G(X) \le \left\{
        \begin{tabular}{ll}
                      ${\sf tp}_{\sone(\opkfd,\op)}(X)$        & \mbox{ if ${\sf tp}_{\sone(\opkfd,\op)}(X)$ is a limit ordinal}\\
                      ${\sf tp}_{\sone(\opkfd,\op)}(X)+\omega$ & \mbox{ otherwise.}\\
        \end{tabular}
\right.
\]

\end{corollary}
We now demonstrate a few features of these results.
\begin{proposition}\label{compactcoverrefinement} Let $X$ and $Y$ be compact metrizable spaces. For each open cover $\mathcal{U}$ of $X\times Y$ there are open covers $\mathcal{A}$ of $X$ and $\mathcal{B}$ of $Y$ such that $\{A\times B:A\in\mathcal{A}\mbox{ and }B\in\mathcal{B}\}$ refines $\mathcal{U}$.
\end{proposition}
{\flushleft{\bf Proof:}} 
Let $d_X$ be a compatible metric on $X$ and let $d_Y$ be a compatible metric on $Y$. Define metric $d$ on $X\times Y$ so that $d((x_1,y_1),(x_2,y_2)) = \sqrt{d_X(x_1,x_2)^2+d_Y(y_1,y_2)^2}$. Then $d$ is a compatible metric on $X\times Y$. 

Consider an open cover $\mathcal{U}$ of $X\times Y$. Since $X\times Y$ is compact choose by the Lebesgue Covering Lemma a positive real number $\delta$ such that for each subset $S\subset X\times Y$ with $d$-diameter less than $\delta$ there is a $U\in\mathcal{U}$ with $S\subset U$. Let $\mathcal{V}$ be a finite open cover of $X$ consisting of sets of $d_X$-diameter less than $\frac{\delta}{2}$, and let $\mathcal{B}$ be a finite open cover of $Y$ consisting of open sets of $d_Y$-diameter less than $\frac{\delta}{2}$. Then $\{A\times B:A\in\mathcal{A}\mbox{ and }B\in\mathcal{B}\}$ is an open cover of $X\times Y$ and consists of sets of $d$-diameter less than $\delta$. Thus, this open cover refines $\mathcal{U}$.
$\square$\\

Since $\polspace$ is not countable dimensional, Lemma \ref{ctbldimgc} implies that $dim_{G}(\polspace)\geq\omega+1$. 

\begin {lemma}
For each $n$, $dim_{G}(\polspace\times \closedunitinterval^n) \leq \omega+n+1$.
\end{lemma}
{\flushleft{\bf Proof:}} 
Write $\closedunitinterval^n = A_0\cup A_1\cup\cdots \cup A_n$ where each $A_i$ is zero-dimensional. 
By Proposition \ref{compactcoverrefinement} we may assume the open covers $\mathcal{U}$ of $\polspace\times\closedunitinterval^n$played by ONE are finite, of the form $\{A\times B:A\in\mathcal{A}\mbox{ and }B\in\mathcal{B}\}$ where $\mathcal{A}$ is an open cover of $\polspace$ and $\mathcal{B}$ is an open cover of $I^n$.  

Now consider the following strategy for TWO: Let $F$ be TWO's winning strategy in $\Gc^{\omega+1}(\op,\op)$ on $\polspace$.
In even indexed innings (including index 0), when TWO considers the moves $(O_0,\cdots,O_{2\cdot n})$ of ONE, TWO plays $T^{2\cdot n}_{\polspace} = F(O^0_{\polspace},O^2_{\polspace},\cdots,O^{2\cdot n}_{\polspace})$ on $\polspace$
Since $A_0\subset\closedunitinterval$ is zerodimensional TWO chooses an open refinement $T^{2\cdot n}_{\closedunitinterval}$ of $O^{2\cdot n}_{\closedunitinterval}$ which covers $A_0$.  TWO plays 
\[
  T_{2\cdot n} = G(O_0, O_1,\cdots,O_{2\cdot n - 1},O_{2\cdot n}) = \{U\times V: U\in T^{2\cdot n}_{\polspace}\mbox{ and }V\in T^{2\cdot n}_{\closedunitinterval}\}.
\]

In odd innings the same plan is used on $A_1$ instead of $A_0$. 

Consider the status after $\omega$ innings have elapsed. The uncovered part is contained in a set of the form $C\times \closedunitinterval^n$ where $C$ is compact and zero-dimensional. Thus it takes at most $n+1$ more innings. This shows that 
$dim_{G}(\polspace\times\closedunitinterval^n)\le \omega+n+1$. 
$\square$

Using Lemma \ref{ctbldimgc} for the Smirnov compactum $\smirnovspace$ we have that ${\sf tp}_{\sone(\opkfd,\op)}(\smirnovspace)=2 < \dim_G(\smirnovspace) = \omega$. Thus the second alternative of Corollary \ref{scvsopkfd} cannot be improved. 

For the Pol compactum $\polspace$ we have ${\sf tp}_{\sone(\opkfd,\op)}(\polspace)=\omega+1$. Thus, using $n=0$, $\dim_G(\polspace)=\omega+1$. By the second alternative of Theorem \ref{scvsopfd} we conclude that $\dim_G(\polspace)\le (\omega+1)+\omega = \omega\cdot 2$. Thus the second alternative of Theorem \ref{scvsopfd} does not give optimal information in all cases. 
For $n>1$ we also have that ${\sf tp}_{\sone(\opkfd,\op)}(\polspace^n) = \omega\cdot n+1$ so that Corollary \ref{scvsopkfd} implies that $\dim_G(\polspace^n) \le \omega\cdot(n+1)$. 

\begin{conjecture}
$\dim_G(\polspace^n) = \omega\cdot n+1$.
\end{conjecture}

Let $X$ be $\oplus_{n=1}^{\infty}\polspace^n$, the topological sum of the finite powers of the Pol compactum. This is a locally compact space. Let $\polspace_\omega = (\oplus_{n=1}^{\infty}\polspace^n)\cup\{p_{\omega}\}$ be the one-point compactification of $X$. Then we have $\dim_G(X)=\omega^2={\sf tp}_{\sone(\opkfd,\op)}(X)$ and $\dim_G(\polspace_{\omega})=\omega^2= {\sf tp}_{\sone(\opkfd,\op)}(\polspace_{\omega})$.

\vspace{0.2in}
Continuing investigating game dimension's product theory we also note:
\begin{theorem}\label{cantorfactor} If $X$ is a compact metrizable space and $\mathbb{C}$ is the Cantor set, then $\dim_G(X) = \dim_G(X\times\mathbb{C})$
\end{theorem}
{\flushleft{\bf Proof:}}
When $\mathcal{U}$ is an open cover of $X\times\mathbb{C}$ we may assume by Proposition \ref{compactcoverrefinement} that $\mathcal{U}$ is finite and that there are finite open covers $\mathcal{U}_X$ of $X$ and $\mathcal{U}_{\mathbb{C}}$ of the Cantor set such that $\mathcal{U} = \{U\times V:U\in\mathcal{U}_X\mbox{ and }V\in\mathcal{U}_{\mathbb{C}}\}$. Since $\mathbb{C}$ is zero dimensional we may assume that $\mathcal{U}_{\mathbb{C}}$ is pairwise disjoint.

Let $\alpha$ denote $\dim_G(X)$. We may assume that $\alpha>\omega$. Let $F$ be TWO's winning strategy in $\Gc^{\alpha}(\op,\op)$ on $X$. 

We now define a winning strategy $G$ for TWO in the game $\Gc^{\alpha}(\op,\op)$ played on $X\times\mathbb{C}$. 
With $\gamma<\alpha$ and open covers $(\mathcal{U}^0,\cdots,\mathcal{U}^{\gamma})$ of $X\times\mathbb{C}$ given, each finite, compute $\mathcal{V}^{\gamma}_{X} = F(\mathcal{U}^0_{X},\cdots,\mathcal{U}^{\gamma}_{X})$, and compute $\mathcal{U}^{\gamma}_{\mathbb{C}}$, and then define
\[
  G(\mathcal{U}^0,\cdots,\mathcal{U}^{\gamma}) = \{U\times V: U\in\mathcal{V}^{\gamma}_{X}\mbox{ and }V\in \mathcal{U}^{\gamma}_{\mathbb{C}}\}.
\]
This defines $G$ for inning $\gamma$ for each $\gamma<\alpha$.

To see that $G$ is a winning strategy, consider a $G$-play of length $\alpha$, ONE's moves denoted $\mathcal{U}^{\gamma}$, $\gamma<\alpha$, and TWO's, $T^{\gamma}$.

Consider an $(x,c)\in X\times\mathbb{C}$. From the definition of $G$ we have an $F$-play
\[
  \mathcal{U}^{\gamma}_{X},\, T^{\gamma}_{X}=F(\mathcal{U}^{\delta}_{X}:\delta\le \gamma),\, \gamma<\alpha
\]
on $X$, won by TWO. Choose $\gamma<\alpha$ with $x\in\bigcup T^{\gamma}_{X}$. Pick a $V\in T^{\gamma}_{X}$ with $x\in V$. Also, choose a $C\in \mathcal{U}^{\gamma}_{\mathbb{C}}$ with $c\in C$. Then $V\times C$ is in $T^{\gamma}$. 

Thus $G$ is a winning strategy for TWO. This shows that $\dim_G(X\times\mathbb{C})\le \alpha$. Since $X$ is homeomorphic to a closed subspace of $X\times\mathbb{C}$, Lemma \ref{subspace} implies that $\alpha = \dim_G(X)\le \dim_G(X\times\mathbb{C})$.
$\square$\\

\begin{theorem}\label{cantorfactor} Let $X$ and $Y$ be compact metrizable spaces. If $\dim_G(X)$ is a limit ordinal and $Y$ is countable dimensional, then $\dim_G(X) = \dim_G(X\times Y)$
\end{theorem}
{\flushleft{\bf Proof:}} Put $\alpha = \dim_G(X)$. Write $\alpha = \bigcup_{n<\infty}S_n$ where each $S_n$ has order type $\alpha$, and for $m<n$, $S_m\cap S_n = \emptyset$. For each $n$ enumerate $S_n$ in an order preserving way as $(s^n_{\gamma}:\gamma<\alpha)$. Write $Y=\cup_{n<\infty}Y_n$ where each $Y_n$ is a zerodimensional set. 

Let $F$ be a winning strategy for ONE in the game $\Gc^{\alpha}(\op,\op)$ on $X$. We define a winning strategy $G$ for $Y$ in the game $\Gc^{\alpha}(\op,\op)$ on $X\times Y$:
First, when ONE plays an open cover $\mathcal{U}$ of $X\times Y$, TWO will replace it, using Proposition \ref{compactcoverrefinement}, with a finite refinement of the form $\{U\times V:U\in\mathcal{U}_X\mbox{ and }V\in\mathcal{U}_Y\}$ where $\mathcal{U}_X$ is an open cover of $X$ and $\mathcal{U}_Y$ is an open cover of $Y$.

Now we define $G$. Let an inning $\gamma<\alpha$ be given, as well as a move $O_{\gamma}$ of ONE. Choose $n$ with $\gamma\in S_n$, say $\gamma = s^n_{\nu}$. Replace $O_{\gamma}$ with the refinement $\{U\times V:U\in O^{\gamma}_X\mbox{ and }V\in O^{\gamma}_Y\}$ where $O^{\gamma}_X$ is a finite open cover of $X$ and $O^{\gamma}_Y$ a finite open cover of $Y$. Let $\mathcal{V}_{\gamma}$ be a finite disjoint refinement of $O^{\gamma}_Y$ which covers $Y_n$, and let $T^{\gamma}_X = F(O^{\gamma}_{X,s^n_0},\cdots,O^{\gamma}_{X,s^n_{\nu}})$. Then define
\[
  T_{\gamma} = G(O_{\nu}:\nu\le\gamma):=\{U\times V:U\in T^{\gamma}_X\mbox{ and }V\in \mathcal{V}_{\gamma}\}.
\]
It follows that $G$ is a winning strategy for TWO in $\Gc^{\alpha}(\op,\op)$, showing $\dim_G(X\times Y)\le \alpha$. But $X$ is homeomorphic to a closed subspace of $X\times Y$, implying by Lemma \ref{subspace} that $\alpha \le \dim_G(X\times Y)$. $\square$\\


\begin{thebibliography}{}
\bibitem{addisgresham} D.F. Addis and J.H. Gresham, \emph{A class of infinite dimensional spaces. Part I: Dimension theory and Alexandroff's problem}, {\bf Fundamenta Mathematicae} 101 (3) (1978), 195 - 205. 
\bibitem{lbscgame} L. Babinkostova, \emph{Selective screenability game and covering dimension}, {\bf Topology Proceedings} 29:1 (2005), 13 - 17.
\bibitem{lbms} L. Babinkostova and M. Scheepers, \emph{Selection Principles and Countable Dimension}, {\bf Note di Mathematica} 27:1 (2007), 5 - 15.

\bibitem{borstcspaces} P. Borst, \emph{Some remarks concerning C-spaces}, {\bf Topology and its Applications} 154 (2007), 665 - 674.

\bibitem{baldwin} S. Baldwin, {\emph Possible Point-Open Types of Subsets of the Reals}, {\bf Topology and its Applications},  Vol.38 NO.3 (1991), 219 - 223.

\bibitem{bernerjuhasz} Berner and I. Juhasz, {\emph  Point-picking games and HFD's}, {\bf Lecture Notes in Mathematics}, Vol. 1103, Springer-Verlag, Berlin-New York (1984), 53-66.


\bibitem{danielsgruenhage} P. Daniels and G. Gruenhage,{\emph  The point-open game on subsets of the reals}, {\bf Topology Appl. 37 (1990)}, 53-64.


\bibitem{engelkingpol} R. Engelking and E. Pol, \emph{Countable dimensional spaces: a survey}, {\bf Dissertationes Mathematicae} 216, (1983).

\bibitem {nagami} K. Nagami, {\emph Monotone Sequence of O-dimensional Subsets of Metric Spaces }, {\bf Proceedings of the Japan Academy}, Vol.41 , No.9(1965) 771-772. 
 
\bibitem{rpol} R. Pol, \emph{On classification of weakly infinite dimensional compacta}, {\bf Fundamenta Mathematicae} 116:3 (1983), 169 - 188.

\bibitem {smirnov} Y. M. Smirnov, {\emph On universal spaces for certain classes of infinite dimensional spaces (English translation)}, {\bf American Mathematical Society Translations Series} 2 21 (1962), 35 - 50.

\bibitem{diaggamelength} M. Scheepers, \emph{The length of some diagonalization games}, {\bf Archive for Mathematical Logic} 38 (1999), 103 - 122. 

\bibitem{succtype} M. Scheepers, {\emph A topological space could have infinite successor point-open types}, {\bf Topology and its Applications}, (1995), 95 - 99.

\end{thebibliography}
\end{document}